\documentstyle[12pt,leqno]{article}
\title{
 \bf Contact topology and $CR$ geometry
in three dimensions}
\author{
  Jih-Hsin Cheng\\
Institute of Mathematics, Academia Sinica\\
Taipei, R.O.C.(Taiwan)\\
E-mail: cheng@math.sinica.edu.tw}
\date{}

\begin{document}
\maketitle
{\normalsize \setlength{\baselineskip}{18pt}
\begin{abstract}
We study low-dimensional problems in topology and geometry via a study
of contact and Cauchy-Riemann ($CR$) structures. A contact structure is called
spherical if it admits a compatible spherical $CR$ structure. We will talk
about spherical contact structures and our analytic tool, an evolution
equation of $CR$ structures. We argue that solving such an equation for the
standard contact 3-sphere is related to the Smale
conjecture in 3-topology. Furthermore, we propose a contact analogue
of Ray-Singer's analytic torsion. This ''contact torsion'' is expected to be
able to distinguish among ''spherical space forms'' $\{\Gamma\backslash
S^{3}\}$ as contact manifolds. We also propose the study of a
certain kind of monopole equation associated with a contact structure. In
view of the recently developed theory of contact homology algebras, we will
discuss its overall impact on our study.
\end{abstract}

\section{Spherical contact structures}

\medskip

Let $(M^{3},\xi)$ denote a contact 3-manifold with the contact structure $\xi
$. (assume $M^{3}$ oriented and $\xi$ cooriented if necessary) We call an
almost complex structure $J$ on $\xi$ a $CR$ (stands for Cauchy-Riemann)
structure (compatible with $\xi$). That is to say, an endomorphism $J:\xi$
$\rightarrow\xi$ \ such that $J^{2}=-Id_{\xi}$. There are no local invariants
for $(M^{3},\xi)$ according to a well known theorem of Darboux. Also for
closed $M^{3}$, two nearby contact structures are isotopy-equivalent by a
theorem of Gray. ([Gr], [Ham]) Therefore a contact structure on a closed
$M^{3}$ has no continuous moduli. On the other hand, we do have local
invariants for a $CR$ 3-manifold $(M^{3},\xi,J)$. Namely, we can talk about
''curvature''. Our strategy of studying 3-topology is via a study of contact
topology and $CR$ geometry.

\medskip

There are distinguished $CR$ structures $J$, called spherical, if
$(M^{3},\xi,J)$ is locally $CR$ equivalent to the standard 3-sphere
$(S^{3},\hat{\xi},\hat{J})$, or equivalently if there are contact coordinate
maps into open sets of $(S^{3},\hat{\xi})$ so that the transition contact maps
can be extended to holomorphic transformations of open sets in $C^{2}$. In
1930's, Elie Cartan ([Ca], [CL1]) obtained a geometric quantity, denoted as
$Q_{J}$, by solving the local equivalence problem for the $CR$ structure so
that the vanishing of $Q_{J}$ characterizes $J$ to be spherical. We will call
$Q_{J}$ the Cartan (curvature) tensor. A contact structure $\xi$ is called
spherical if there is a spherical $CR$ structure compatible with it.

\medskip

Our main concern is the existence problem of spherical contact
structures. For instance, we ask if any homology 3-sphere admits a spherical
contact structure, or does there exist a nonspherical contact homology
3-sphere? Notice that a spherical (contact) homotopy 3-sphere is
contact-diffeomorphic to $(S^{3},\hat{\xi})$.

\medskip

It has been believed that for closed $M^{3}$, a spherical (contact)
structure is tight. (for open $M^{3}$, Eliashberg gave counterexamples) We
probably can prove this conjecture by showing that the contact homology
(recently developed by Eliashberg, Givental and Hofer, [EGH]) of a spherical
structure does not vanish. So, most likely, $\{$spherical structures$\}$ is a
restricted class of tight contact structures, which we can apply more
analytic tools to study.

\bigskip

\section{The Cartan flow}
\setcounter{equation}{0}

\medskip

The tool we'd like to use is the so-called Cartan flow, an evolution equation
for $CR$ structures $J_{(t)}$ on $(M^{3},\xi)$:

\begin{equation}
\partial_{t}J_{(t)}=Q_{J_{(t)}}.
\end{equation}

\noindent
Namely, we deform a $CR$ structure in the direction of its Cartan tensor. And
we hope that the limit $CR$ structure has the vanishing Cartan tensor,
therefore is spherical. First note that (2.1) is a system of 4-th order
nonlinear subparabolic equations (up to an action of contact diffeomorphisms).
([CL1]) Second, we'll mention some topological and geometrical implications of
solving such an equation. Before doing that, let's see what have been known.

\medskip

In the late 1980's, it was observed that (2.1) is a downward
(negative) gradient flow. In fact, D. Burns and C. Epstein (also J. Lee and
myself) ([BE], [CL1]) found an energy functional $\mu_{\xi}$ defined on a
certain space of $CR$ structures (assuming trivial holomorphic tangent bundle for instance) so that

$$\delta\mu_{\xi}(J)=-Q_{J}$$

\noindent
(meaning $D\mu_{\xi}(J)(E)=-<Q_{J},E>$ for any tangent vector $E$ at $J$, in
which $<$ $,>$ is the inner product induced by the Levi form).

\medskip

The short time solution can be proved by adding a gauge-fixing term
to the right-hand side of (2.1). The linearization of the resulting equation
is subparabolic with the leading space term of the form $-(const)L_{\alpha
}^{\ast}L_{\alpha}u.$Here $L_{\alpha}$ is the generized Folland-Stein operator
and subelliptic if $\alpha$ is not an odd integer. In our case, $\alpha
=4+i\sqrt{3}$. ([CL1])

\medskip

Now we come back to the first potential application in 3-topology by
solving (2.1) just for $(M^{3},\xi)=(S^{3},\hat{\xi})$. This will confirm the
so-called Smale conjecture: $Diff(S^{3})\approx O(4)$ ($"\approx"$ means
''homotopy equivalent'') as first pointed out by Eliashberg in the early
1990's. In fact, Hatcher ([Hat]) gave a combinatorial proof in 1983. But
people are still seeking for more geometric proofs. We can argue that the
solution for $(S^{3},\hat{\xi})$ implies the Smale conjecture as follows.
Since $\hat{J}$ is the unique spherical $CR$ structure on $(S^{3},\hat{\xi})$,
any other $J$ will converge to $\hat{J}$ through the Cartan flow. This means a
certain marked $CR$ moduli space $\Im^{^{\prime}}/C^{\prime}$ is contractible.
But $\Im^{^{\prime}}$, the space of certain marked $CR$ structures, is
contractible. So $C^{\prime}$, the group of certain marked contact
diffeomorphisms, is contractible too. It follows that $Diff(S^{3})\approx
O(4)$ by the relation between $Diff(S^{3})$ and $C^{\prime}$.

\medskip

We remark that $\hat{J}$ is a strict local minimum for $\mu_{\hat{\xi
}}$. ([CL2]) The solution to (2.1) for $(S^{3},\hat{\xi})$ will imply that
$\hat{J}$ is actually a global minimum. Because there seems to be no suitable
maximum principle available for 4-th order subelliptic operators, a proof of
the solution to (2.1) would probably have to be based on a priori integral
estimates in place of the usual pointwise estimates for 2nd order parabolic
flows. To learn more analytic techniques, we have been working on a
comparatively easier flow. Let us define an energy $e_{J}$ for a contact form
$\theta$ as follows:

$$
e_{J}(\theta)=\int_{M^{3}}(W_{J,\theta})^{2}%
\theta\wedge d\theta.
$$

\noindent
Here $W_{J,\theta}$ denotes the Tanaka-Webster curvature associated with
$(J,\theta)$. ([Ta], [We]) We consider the downward gradient flow of $e_{J}$.
If we write $\theta_{(t)}=e^{2\lambda_{(t)}}\hat{\theta}$ with respect to a
fixed background contact form $\hat{\theta}$, then the equation can be
expressed in $\lambda_{(t)}$ as

\begin{equation}
\partial_{t}\lambda_{(t)}=\Delta
_{b}W_{J,\theta_{(t)}}.
\end{equation}

\noindent
Here $\Delta_{b}$ denotes the (positive) sublaplacian. (notice the sign
difference for $\Delta_{b}$ in [Lee]) The equation (2.2) is a 4-th order
subparabolic, but scalar, flow. (while (2.1) is a ''vector'' flow with two
independent real unknowns) It is easy to see that the volume $\int_{M^{3}%
}\theta_{(t)}\wedge d\theta_{(t)}=\int_{M^{3}}e^{4\lambda_{(t)}}\hat{\theta
}\wedge d\hat{\theta}$ is preserved under the flow (2.2). Under certain
conditions, we can establish the following integral estimate: ([CCg])

\begin{equation}
\partial_{t}\int_{M^{3}}%
e^{5\lambda_{(t)}}\hat{\theta}\wedge d\hat{\theta}\leq C.
\end{equation}

\noindent
Here the constant $C$ may or may not depend on the maximum time according to
the applied situations. 
%The $C^{0}$ estimate, hence the long time solution and asymptotic convergence, follows from (2.3) and other known techniques.
The idea of estimating an integral such as the one in (2.3) comes from
the study of a certain metric flow related to general relativity. 
The involved integral
quantity is known as the Bondi mass. We wonder if there are Bondi-mass type
estimates for the Cartan flow (2.1).

\medskip

A $CR$ manifold is embeddable if it can be ''realized'' as the
boundary of a compact complex manifold. (with the $CR$ structure being the one
induced from the complex structure) The embeddability is a special property
for 3-dimensional $CR$ manifolds since any closed $CR$ manifold of dimension
$\geq$ $5$ is embeddable. ([BdM]) Now it is natural to ask the following question:

\medskip

{\bf{Is the embeddability preserved under the Cartan flow (2.1)?}}

\medskip

By a direct construction of an integrable almost complex structure, we
can show that if $J_{(0)}$ is embeddable with the torsion $\equiv$
$L_{T}J_{(0)}=0$ and $W_{J_{(0)},\theta}>0$ (or $<0$), then $J_{(t)}$ stays
embeddable (for a short time). ([Ch2]) Here $T$ denotes the Reeb vector field
associated with $\theta.$ In fact, the torsion stays zero under the flow. Also
the existence of a $CR$ vector field $T$ is sufficient to imply the
embeddability of the $CR$ structure as pointed out by L\'{a}szl\'{o} Lempert. ([Lem])
So the condition on the Tanaka-Webster curvature is redundant. We conjecture that the
embeddability is preserved under the Cartan flow without any conditions.

\medskip

On the other hand, the zero torsion condition reduces the complexity
of our flow a lot. It seems to be a good starting point. We are in a situation
analogous to Hamilton's Ricci flow. Namely given a closed contact 3-manifold
$(M^{3},\xi)$. Suppose there is a pseudohermitian structure $(J,\theta)$ with
vanishing torsion and positive Tanaka-Webster curvature. Then can we conclude
that $\xi$ is spherical? A possible proof is to apply the Cartan flow to show
that the limit $CR$ structure (together with the fixed contact form $\theta$)
has the positive constant Tanaka-Webster curvature. (recall that the torsion
stays zero for all time) Therefore it has the vanishing Cartan tensor. So it
is spherical.

\bigskip

\section{Spherical space forms}

\medskip

Since the linearization $\delta Q_{J}$ of the Cartan tensor $Q_{J}$ is
subelliptic modulo the action of the contact diffeomorphism group $C_{\xi}$
([CL1]), the kernel of $\delta Q_{J}$ is finite-dimensional modulo the action
of $C_{\xi}$. So the ''virtual'' dimension of the moduli space of spherical
$CR$ structures is finite-dimensional. In this section, we will just consider
a class of examples for the 0-dimensional case. Let $\Gamma$ denote a fixed
point free finite subgroup of the $CR$ automorphism group of the standard
$S^{3}$(which is isomorphic to $PU(2,1)$). Then the quotient space
$\Gamma\backslash S^{3}$ inherits a (spherical) contact structure from
$(S^{3},\hat{\xi}).$ It's natural to work on the following problem.

\medskip

{\bf{Problem: Classify $\{\Gamma\backslash S^{3}\}$ as contact manifolds.}}

\medskip

It has been believed that $\Gamma_{1}\backslash S^{3}$ is
contact-diffeomorphic to $\Gamma_{2}\backslash S^{3}$ if and only if they are
$CR$-diffeomorphic to each other in analogy with the conformal case. Thus to
deal with the above problem, we borrow ideas from quantum physics to find a
potential invariant in terms of $CR$ geometry.

\medskip

If we view $\mu_{\xi}$ as a Lagrangian (action, more accurately) in
$2+1$ dimensions, spherical $CR$-structures are just classical fields.
Therefore, ``quantum fluctuations'' should give us refined invariants.\ In
practice, we compute the partition function heuristically:

\begin{eqnarray}
{\begin{array}{rl}
{\cal{Z}}_{k}  & ={\int_{\Im_{\xi}/\cal{C}_{\xi}}}{\cal{D}}
[J]e^{ik{\mu_{\xi}}([J])}\nonumber\\
& =k^{-\frac{dim}{2}}({\cal{Z}}_{sc}+O(k^{-1}))\:(k\:large),\nonumber
\end{array}
}
\end{eqnarray}

\noindent in which ${\cal{Z}}_{sc}$ is called the semi-classical approximation. Note
that only classical fields make contributions to ${\cal{Z}}_{sc}$. By imitating the
finite dimensional case, we can compute the modulus of ${\cal{Z}}_{sc}$:

\begin{eqnarray}
{
\begin{array}{rl}
|{\cal{Z}}_{sc}|  & =lim_{k{\rightarrow}{\infty}}k^{\frac{dim}{2}%
}|{\cal{Z}}_{k}|\nonumber\\
& ={\Sigma}_{J:spherical}\left|  \frac{det{\Box_{J}}}{det^{\prime}{\delta
}Q_{J}}\right|  ^{\frac{1}{2}},\nonumber
\end{array}
}
\end{eqnarray}

\noindent in which $\Box_{J}$ is a fourth-order subelliptic self-adjoint
operator related to the $C_{\xi}$-action, and $\delta Q_{J}$, the second
variation of $\mu_{\xi}$, is also a fourth-order subelliptic self-adjoint
operator modulo the $C_{\xi}$-action. We can regularize two determinants via
zeta functions. ($det^{\prime}$ means taking a regularized determinant under a
certain gauge-fixing condition.) (see [Ch1] for more details)

\medskip

{\bf{\noindent Conjecture: If $J$ is spherical,}}

$$
Tor(J){\equiv}\left|  \frac{det{\Box_{J}}}{det^{\prime}{\delta}
Q_{J}}\right|  ^{\frac{1}{2}}
$$

{\bf{\noindent is independent of any choice of contact form, i.e., a $CR$ invariant.}}

\medskip

We expect to use $Tor(J)$ to distinguish among spherical space forms
$\{\Gamma\backslash S^{3}\}$. And we note that $Tor(J)$ is a contact-analogue
of Ray-Singer's analytic torsion while no contact-analogue is known for the
Reidemeister torsion. Also we speculate that if the contact homology of
$\Gamma\backslash S^{3}$ ([EGH]) can distinguish $\Gamma\backslash S^{3}$'s,
it may be possible to identify $Tor(J)$ with a certain quantity composed of
elements in the contact homology of $\Gamma\backslash S^{3}$.

\medskip

Let us consider the case that $\Gamma=I^{\ast}$, the binary icosahedral
group. It is known that $I^{\ast}\backslash S^{3}$ is just the Poincar\'{e}
homology sphere $P$. Therefore its contact structure is spherical. We know
that two spherical $CR$ manifolds can be glued together to form the spherical
connected sum by an orientation-preserving gluing map. (the gluing map is
given by a $CR$ inversion defined on the Heisenberg group $H$ minus the origin
in view of a spherical $CR$ manifold being $CR$ equivalent to $H$ locally. In
coordinates $(t,z)$ where $t\in R,z\in C$, the $CR$ inversion $I$ defined by
$I(t,z)=(-t/|w|^{2},z/w)$ in which $w=t+i|z|^{2}$ satisfies $I^{\ast}%
\theta_{0}=|w|^{-2}\theta_{0}$ where $\theta_{0}=dt+izd\bar{z}-i\bar{z}dz$ is
the standard contact form on $H$. It is easy to verify that ''$I$'' is
orientation preserving and interchanges the surfaces defined by $|w|=2$ and
$|w|=1/2$, respectively). It follows that the connected sum $P\#P$ of $P$ and
itself is spherical too. On the other hand, we have the following conjecture

\medskip

{\bf{Conjecture: There does not exist any spherical contact structure
on $P\#\bar{P}$.}}

\medskip

\noindent
Here $\bar{P}$ denotes $P$ with the reverse orientation. In [EH], Etnyre and
Honda proved that there does not exist any tight contact structure on
$P\#\bar{P},$ either positive or negative. So if a spherical contact structure
on a closed 3-manifold is tight (a previous conjecture that we mentioned in
the end of section 1), then the above-mentioned conjecture holds in view of
Etnyre and Honda's result. If so, we then have a homology 3-sphere that does
not admit any spherical contact structure.

\bigskip

\section{Monopoles and contact structures}

\medskip

Given a contact 3-manifold $(M^{3},\xi)$ and a background pseudohermitian
structure $(J,\theta)$, we can discuss a canonical $spin^{c}$-structure
$c_{\xi}$ on $\xi^{\star}$. ([CCu]) With respect to $c_{\xi}$, we will
consider the equations for our ``monopole'' $\Phi$ coupled to the ``gauge
field'' $A$. Here, $A$, the $spin^{c}$-connection, is required to be
compatible with the pseudohermitian connection on $M^{3}$. The Dirac operator
$D_{\xi}$ relative to $A$ is identified with a certain boundary $\bar
{\partial}$-operator $\sqrt{2}(\bar{\partial}_{b}^{a}+(\bar{\partial}_{b}%
^{a})^{\star})$. In terms of the components $(\alpha,\beta)$ of $\Phi$, our
equations read as

\medskip

\begin{flushleft}(4.1)\hspace{.8in}
$\left\{
\begin{array}[c]{c}
({\bar{\partial}}_{b}^{a}+({\bar{\partial}}_{b}^{a})^{\star})({\alpha}+{\beta
})=0\\
(or\:{\alpha}_{,{\bar{1}}}^{a}=0,\:{\beta}_{{\bar{1}},1}^{a}=0)\\
da(e_{1},e_{2})-W_{J,\theta}=|{\alpha}|^{2}-|{\beta}_{\bar{1}}|^{2},
\end{array}
\right.  $
\end{flushleft}

\medskip

\noindent where $A=A_{can}+iaI$ and $W_{J,\theta}$ denotes the Tanaka-Webster
curvature. Our first step in understanding (4.1) is as follows:

\medskip

{\it Suppose the torsion $L_{T}J=0$ ($T$ is the
Reeb vector field). Also, suppose $\xi$ is symplectically
semifillable, and that the Euler class $e(\xi)$ is not a torsion
class. Then (4.1) has nontrivial solutions (i.e., $\alpha$ and
$\beta$ are not identically zero simultaneously). ([CCu])}

\medskip

On the other hand, the Weitzenbock-type formula gives a nonexistence
result for $W_{J,\theta}>0$. Together with the above existence result, we can
conclude the following:

\medskip

{\it Suppose the torsion vanishes and the Tanaka-Webster curvature
$W_{J,\theta}>0$. Then, either $\xi$ is not symplectically
semifillable, or $e(\xi)$ is a torsion class. ([CCu])}

\medskip

We remark that Rumin ([Ru]) proved that $M^{3}$ must be a rational
homology sphere under the conditions given above using a different method.
Originally we were hoping to define contact invariants from the solution space
of (4.1). But since $D_{\xi}$ (also $da(e_{1},e_{2})$) is not elliptic (not
even subelliptic), the solution space might be infinite dimensional. To
distinguish such spaces, it seems that we need to know more structures about
the solution space. On the other hand, the contact homology algebras recently
developed in [EGH] seem to provide such a structure from the algebraic point
of view.

\bigskip

\section{General discussion}

\medskip

About the Cartan flow (2.1), one would like to know under what conditions the
solution to (2.1) exists for all time and converges as $t\rightarrow+\infty$
to a spherical $CR$ structure. This will be impossible in general. Even if our
manifold is the sphere, if we start with an overtwisted contact structure, the
solution to (2.1) can not converge since the limit $CR$ structure would
perforce be diffeomorphic to the standard one (which is tight). Hence the
solution must blow up at a finite time. We then ask what the shape of the
blow-up set looks like.

\medskip

In [Go], W. Goldman obtained some topological obstruction to the
existence of spherical (contact) structures. In particular, $T^{3}$ does not
admit any spherical structures. We hope to be able to obtain some contact
topological obstruction in terms of contact homology algebras. To do this, we
have to analyze how the contact homology changes under covering and developing
maps associated with a spherical structure. K. Mohnke [Mo] has studied the
contact homology of certain coverings. His work should be useful for our
study. Haven't obtained some contact obstruction, we can then answer the
nonexistence problem of spherical structures in a more refined way. For
instance, we might be able to determine which ones among those known tight
contact structures are nonspherical for Brieskorn homology spheres
$\sum(2,3,6n-1),n\geq2$. \ Also we can then easily confirm the following
previously mentioned conjecture by showing that the contact homology of a
spherical structure does not vanish.

\medskip

{\bf Conjecture: A spherical structure on a closed $3$-manifold is tight.}

\bigskip

{\bf{Acknowledgments.}} In preparing this article the author benefited from a number of
conversations with Yasha Eliashberg, L\'{a}szl\'{o} Lempert, Kaoru Ono, Klaus
Mohnke, and Mei-Lin Yau. Also, the author would like to thank the organizers
of both the 3rd Asian Mathematical Conference (Quezon City) and the Contact
Geometry Conference (Stanford) for their kind invitations. The main part of
this article is based on the author's notes for his talks at the above
mentioned two conferences.The research was partially supported by the National
Science Council, R.O.C., under grant NSC 90-2115-M-001-003.

\end{document}